\documentclass[12pt]{amsart}
\usepackage{amssymb, epic,eepic,epsfig,amsbsy,amsmath,amscd}
\usepackage[left=2cm,right=2cm,top=3cm,bottom=3.5cm,a4paper]{geometry}
\usepackage{color}
\usepackage{verbatim} % To be able to use the  which "outcomments" some part of the file.

\usepackage[small,nohug,neads=littlevee]{diagrams}
\diagramstyle[labelstyle=\scriptstyle]

\def\be#1\ee{\begin{equation}#1\end{equation}}
\theoremstyle{plain}

\newtheorem{Thm}{Theorem}

\newtheorem{proposition}[Thm]{Proposition}%[section]
\newtheorem{lemma}[Thm]{Lemma}
\newtheorem*{lemma*}{Lemma}
\newtheorem*{Theorem*}{Theorem}

\newtheorem*{conjecture}{Conjecture}

\theoremstyle{definition}

\newtheorem*{remark*}{Remark}
\newtheorem*{example*}{Example}

\def\printname#1{
    \if\draft
        \smash{\makebox[0pt]{\hspace{-0.5in}
            \raisebox{8pt}{\tt\tiny #1}}}
    \fi
}

\newlength{\standardunitlength}
\catcode`\@=11
\long\def\@makecaption#1#2{%
     \vskip 10pt

\setbox\@tempboxa\hbox{%\ifvoid\tinybox\else\box\tinybox\fi
       \small\sf{\bfcaptionfont #1. }\ignorespaces #2}%
     \ifdim \wd\@tempboxa >\captionwidth {%
         \rightskip=\@captionmargin\leftskip=\@captionmargin
         \unhbox\@tempboxa\par}%
       \else
         \hbox to\hsize{\hfil\box\@tempboxa\hfil}%
     \fi}
\font\bfcaptionfont=cmssbx10 scaled \magstephalf
\newdimen\@captionmargin\@captionmargin=2\parindent
\newdimen\captionwidth\captionwidth=\hsize
\catcode`\@=12

\newcommand{\ord}{\operatorname{ord}}

\def\F{\mathbb F}

\def\BZ{\mathbb Z}

\def\BQ{\mathbb Q}

\def\BC{\mathbb C}

\def\cT{\mathcal T}
\def\cM{\mathcal M}

\def\A{\mathcal A}

\def\cR{\mathcal R}

\def\lk{\operatorname{lk}}

\def\cL{\mathcal L}

\def\M{\mathcal M}
\def\N{{\mathbb N}}
\def\F{\mathcal F}

\def\R{\mathcal R}

\def\calS{\mathcal S}

\newcommand{\sn}{\operatorname{sn}}

\def\bn{{\mathbf n}}

\def\lk{{\text{lk}}}

\def\Z{\BZ}

\def\bk{\mathbf{k}}

\def\Q{\BQ}
\def\n{{\mathbf n}}

\def\ev{\mathrm{ev}}

\def\Habiro{\widehat{\Z[q]}}

\newcommand\nc{\newcommand}

\nc\FIG[3]{\begin{figure}
    \includegraphics[#3]{#1.eps}
    \caption{#2}
    \label{fig:#1}
    \end{figure}}

\nc\incl[2]{{\includegraphics[height=#1]{#2.eps}}}

\nc\zzzcolon {\colon\thinspace}
\nc\zzzvert {\ |\ }
\nc\simh{\underset{h}{\sim}}

\nc\trr{\triangleright}
\nc\et{\incl{1em}{bot0}}
\nc\ul{\underline}

\nc\modone {{\mathbf 1}}
\nc\modA {{\mathcal A}}
\nc\sfA{{\mathsf A}}
\nc\modb {{\mathsf b}}
\nc\modB {{\mathsf B}}
\nc\cB{\check{\modB }}
\nc\modC {{\mathcal C}}
\nc\modH {{\langle \HH\rangle }}
\nc\modh {{\mathsf h}}
\nc\modI {{\mathcal I}}
\nc\modJ {{\mathsf J}}
\nc\modk {{\mathbf k}}
\nc\modL {{\mathcal L}}
\nc\LL{{\mathsf L}}
\nc\HH{{\mathsf H}}
\nc\modM {{\mathcal M}}
\nc\sfP{{\mathsf P}}
\nc\modP {{\mathcal P}}
\nc\modQ {{\mathbb Q}}
\nc\modR {{\mathbb R}}
\nc\modr {v}
\nc\modS {{\mathcal S}}
\nc\modT {{\mathcal T}}
\nc\modV {{\mathsf V}}
\nc\modY {{\mathcal Y}}
\nc\modZ {{\mathbb Z}}
\nc\ad{{\operatorname{ad}}}
\nc\coad{{\operatorname{coad}}}
\nc\coadb{\underline{\coad}}
\nc\adb{\underline{\ad}}
\nc\coev{{\operatorname{coev}}}
\nc\pat{{\operatorname{\sf pat}}}
\nc\low{{\operatorname{low}}}
\nc\up{{\operatorname{up}}}

\nc\evdn{{\incl{.8em}{evdown}}}
\nc\cvdn{{\incl{.8em}{coevdown}}}
\nc\evup{{\incl{.8em}{evup}}}
\nc\cvup{{\incl{.8em}{coevup}}}

\nc\bH{{\underline H}}
\nc\bD{{\underline\Delta }}
\nc\bS{{\underline S}}
\nc\tS{{\tilde{S}}}

\nc\Ob{\operatorname{Ob}}

\newcommand{\psdiag}[3]{\hspace{1mm}\raisebox{-#1mm}{\epsfysize#2mm
\epsffile{#3.eps}}\hspace{1mm}}

\begin{document}

\title[Unified  3--manifold  invariants]{On the
unification of
quantum 3--manifold  invariants}

\author{Anna Beliakova}
%\author{Irmgard B\"uhler}
\address{Institut f\"ur Mathematik, Universit\"at Zurich,
Winterthurerstrasse 190, 8057 Z\"urich, Switzerland}
\email{anna@math.uzh.ch, irmgard.buehler@math.uzh.ch}

\author{Thang Le}
\address{Department of Mathematics, Georgia Institute of Technology,
Atlanta, GA 30332--0160, USA }
\email{letu@math.gatech.edu}

\begin{abstract}
In 2006 Habiro initiated a construction of generating functions for
Witten--Reshetikhin--Turaev (WRT) invariants
known as unified WRT invariants. In a series of papers
together with Irmgard B\"uhler
and Christian Blanchet
we extended his construction to a larger class of 3--manifolds.
The unified invariants provide a strong tool to study  properties
of the whole collection of WRT invariants, e.g.  their integrality,
and hence, their categorification.

In this paper we give   a survey
on ideas and techniques used in the construction of the unified
invariants.

%New structural properties of the
%set of quantum invariants at roots of unity not coprime
%to the torsion are the main applications of our construction.

\end{abstract}

\maketitle

%%%%%%%%%%%%%%%%%%%%%%%%%%%%%%%%%%%%%%%%%%%%%%%%%%%%%%%%%%%%%%%%%%%%%%%%%%%%%%%%%%%%%%%%%%%%%%%%%%%%%%%%%%%%%%%%%%%%%%%%%%%%%%%%%%%%%%%%%%%%%%%%%%%%%%%%%%%%%%%%%%%%%%%%%%%%%%%%%%%%%%%%

\vskip1mm


\section*{Introduction}

\subsection*{Background}
In the 60s and 70s Rochlin, Lickorish and Kirby established a remarkable
connection between links and 3--manifolds.
They showed that every 3--manifold can be obtained by surgery (on $S^3$) along framed links,
and surgeries along two links give the same 3--manifold if and only if the links
are related by a sequence of  Kirby moves. This allows us to think of  3--manifolds as
 equivalence classes of  framed links modulo the relation generated by Kirby moves.

After
the discovery of the Jones polynomial in 1984,  knot theory experienced
the transformation
from an esoteric branch of pure mathematics  to a   dynamic
research field with deep connections to  mathematical physics,
the theory of integrable and dynamic systems, von Neumann algebras,
representation theory, homological algebra, algebraic geometry, etc.
Among important developments were the constructions of {\em quantum link invariants} (generalization of the Jones polynomial), and
of {\em the Kontsevich integral} (universal finite type invariant).

The quantum link invariants were extended to 3--manifolds by
Witten and Reshetikhin--Turaev (WRT), the lift of the
Kontsevich integral to 3--manifolds was defined by Le--Murakami--Ohtsuki (LMO).
However, while the relationship between the Kontsevich integral and  quantum link invariants is simple,
the relationship between the LMO invariant and the WRT invariants  is much more complicated and remains mysterious in many cases.

Let us explain this in more details.
The Kontsevich integral of a knot
takes values in a certain  algebra $\mathcal A(S^1)$ of chord diagrams.
Any semi--simple Lie algebra and  its module define a map $\mathcal A(S^1) \to \Q[[h]]$
called a weight system. One important result of Le--Murakami and Kassel is
that the following diagram commutes, where the $\frak{sl}_2$ weight system  uses the 2--dimensional defining representation.

\begin{diagram}
                &                                         &\BZ[q^{\pm 1}]&                                 &\\
                &\ruTo^{\text{Jones polynomial}}           &              &\rdTo^{q=e^h}                    &\\
\{\text{knots}\}&                                         &              &                                 &\BQ[[h]]\\
                &\rdTo_{\text{Kontsevich integral}}&              &\ruTo_{\frak {sl}_2\text{-weight system}}&\\
                &                                         &\A(S^1)            &                                 &
\end{diagram}

In particular, this proves that the Kontsevich integral dominates the Jones polynomial, and similarly all
quantum link invariants coming from
Lie algebras. Hence in addition to being universal for finite type, the Kontsevich integral is also universal
for all quantum link invariants. It is conjectured that
the Kontsevich integral separates knots.

Does there exist a similar commutative diagram for 3--manifolds?
The  quantum WRT invariant
associates with a compact orientable 3--manifold $M$, a root of unity
$\xi$ and a semi--simple Lie algebra, say $\frak {sl}_2$
 for simplicity,
a complex number $\tau_M(\xi)$. The LMO invariant takes values
in a certain algebra $\mathcal A$ of Feynman diagrams. Every semi--simple Lie algebra defines a weight system map from $\A$ to $\BQ[[h]]$.
Hence, we have

\begin{diagram}
                &                                         &\BC&                                 &\\
                &\ruTo^{\tau_M(\xi)}           &              &\luDotsto                    &\\
\{\text{3--manifolds}\}&                                         &              &                                 &\BQ[[h]]\\
                &\rdTo_{\text{LMO}}&              &\ruTo_{\frak {sl}_2\text{--weight system}}&\\
                &                                         &\A            &                                 &
\end{diagram}

The image of the composition of the two bottom arrows is
known as the Ohtsuki series \cite{Oh}. Ohtsuki showed that for
any rational homology 3--sphere $M$ and a root of unity $\xi$
of {\it prime} order $p$, the first $(p-1)/2$ coefficients, modulo $p$,
of the Ohtsuki series
are determined by $\tau_M(\xi)$. This is
shown by the dotted arrow in the diagram.

This result of Ohtsuki raises many interesting questions: Does there exist
any relationship between the LMO and WRT invariants at roots of
unity of non--prime order? Is the whole set of WRT invariants determined by the LMO invariant?
The discovery of the series prompted Ohtsuki to build the theory of
finite type invariants of homology 3--spheres, which was further developed by
Goussarov and Habiro.

In the 90s, Habegger, Garoufalidis and Beliakova
showed that the LMO invariant is trivial if the first Betti number of a 3--manifold
is bigger than 3, and for the Betti numbers 1,2 and 3,
the LMO invariant is determined
by the classical Casson--Walker--Lescop invariant. However, in the case
of rational homology 3--spheres, the LMO invariant is more powerful, it is 
a universal finite type invariant in the Goussarov-Habiro sense.

The  relationship between the LMO and WRT invariants
at non--prime roots of unity remained open
for quite a while.
This is because all known techniques heavily
rely on the fact that the order of the root is prime and can not
be extended to other roots.  However, recently
 Habiro's theory \cite{Ha} of unified
invariants  provided
a complete solution of this problem in the case
of integral homology 3--spheres. In this case, also
the  question of  integrality of the WRT invariants at non--prime roots of unity
was solved simultaneously. Though
intensively studied (see
\cite{Mu}, \cite{MR}, \cite{GM}, \cite{Le10} and the references there),
 the integrality of WRT invariant was previously known
  for prime roots of unity only.
Note that a conceptual solution of the integrality problem is of primary
importance for any attempt of  categorification of the
WRT invariants (compare \cite{Kho}).

\subsection*{Unified invariants of integral homology 3--spheres}
The unification of the WRT invariants was initiated  in 2006 by Habiro.
For any  integral homology
3--sphere $M$,
Habiro \cite{Ha} constructed a {\em unified invariant} $J_M$
 whose evaluation  at any root of unity coincides with
the value of the  WRT
invariant at that root.
Habiro's unified invariant $J_M$ is an element of the following
ring (Habiro's ring)
\[
\Habiro:=\lim_{\overleftarrow{\hspace{2mm}k\hspace{2mm}}}
\frac{\Z[q]}{
((q;q)_k)},  \qquad \text{ where} \quad (q;q)_k = \prod_{j=1}^k (1-q^j).
\]
Every element $f(q)\in \Habiro$ can be written (non--uniquely) as an infinite sum
\[
f(q)= \sum_{k\ge 0} f_k(q)\, (1-q)(1-q^2)...(1-q^k),
\]
with $f_k(q)\in \Z[q]$. When $q=\xi$, a root of unity,
 only a finite number of terms on the right hand side are not zero,
hence  the evaluation
$\ev_\xi(f(q))$ is well--defined and is an algebraic integer.
However,  the fact that the unified invariant belongs to the Habiro ring
is stronger than just integrality of $\tau_M(\xi)$.

The Habiro ring has beautiful
arithmetic properties.
Every element $f(q) \in \Habiro$ can be considered
 as a function whose domain is the set
 of roots of unity.
%Although each element $f\in \Habiro$ is a function
% defined only at roots of unity,
Moreover, there is a natural Taylor series for $f$ at every root
of unity.
In \cite{Ha} it is shown that two elements $f,g \in \Habiro$ are the same if and
only if their Taylor series at a root of unity coincide.
In addition, each function $f(q) \in \Habiro$ is totally determined
by its values at, say,
infinitely many  roots of order $3^n,\, n\in \N$.
Due to these properties the Habiro ring is also called
a ring of ``analytic functions at roots of unity''. Thus belonging to $\Habiro$ means that the
collection of the  WRT invariants is far from a random
collection of algebraic integers; together they form a nice function.

General properties of the Habiro ring imply that for any integral homology
3--sphere $M$,
the Taylor expansion of the unified invariant $J_M$
at $q= 1$ coincides with the Ohtsuki series
and dominates WRT invariants of $M$ at all roots
of unity (not only of prime order). This is summarized in the following
commutative diagram.

\begin{diagram}
                &                   &\BC             &                                           &\\
                &\ruTo^{\tau_M(\xi)}&\uTo_{q=\xi}    &                                 &\\
\{\BZ\text{HS}\}&\rTo^{J_M(q)}      &\widehat{\BZ[q]}&\rInto                                     &\BZ[[1-q]] \\
                &\rdTo_{\text{LMO}} &                &                                           &\dInto_{h=1-q} \\
				&                   &\A              &\rTo^{\frak {sl}_2\text{--weight}}_{\text{system}} &\BQ[[h]]
\end{diagram}
By $\BZ\text{HS}$ we denoted here the set of integral homology 3--spheres.
In particular, this shows that Ohtsuki series has integral coefficients (which was conjectured by Lawrence and first proved by Rozansky).

Recently, Habiro ring found  applications in algebraic geometry
for  constructing
varieties over the non--existing field of one element \cite{Ma}.

\subsection*{Unified invariants of rational homology 3--spheres}
In \cite{BBL}, we give a full generalization of the Habiro theory to
rational homology 3--spheres. This requires  completely new
techniques coming from number theory, commutative algebra,
quantum group and knot theory.
Let us explain this in more details.

Assume $M$ is a rational homology 3--sphere with
$|H_1(M,\Z)|=b$, where for a finite group $G$ we denote
by $|G|$ the number of its elements.
 Then
our unified invariant $I_M$ belongs to a modification $\R_b$ of the Habiro ring
where $b$ is inverted. Unlike the case  $b=1$, the modified Habiro ring
is not an integral domain, but a direct product of integral domains,
where each factor is determined by its own Taylor
expansion at some root of unity. There is a decomposition
 $I_M=\prod_{c|b} I_{M,c}$, where $I_{M,c}$ dominates the set
$\{\tau_M(\xi)|(\ord(\xi),b)=c\}$. If $b=1$, then $I_M$ coincides with Habiro's $J_M$.
The invariant $I_{M,1}$ was first defined in \cite{Le}.

Our results can be summarized in the following commutative
diagram. Here we assume for simplicity
that $b=p^k$ is a power of a prime and put $e_n:=\exp(2\pi I/n) $
the primitive $n$th root of unity.

\begin{diagram}
                &                        &\BC             &                                          &                 &\\
                &\ruTo(2,4)^{\tau_M(\xi)}&\uTo_{q=\xi}    &                                 &                 &\\
                &                        &\R_b&\rInto                                    &
                                                                                                 \prod^\infty_{i=0}\BZ\left[\frac{1}{p}, e_{p^i}
\right]\left[\left[e_{p^i}-q\right]\right]&\\
                &\ruTo~{I_M(q)}          &                &    &                 &\\
\{\BQ\text{HS}\}&                        &                &                                          & \dTo_{\text{projection to}\; i=0}    &\\
                &\rdTo_{\text{LMO}}      &                &                                          &                 &\\
				&                        &\A              &\rTo^{\frak {sl}_2\text{--weight}}_{\text{system}}&\BQ[[h]]     &
\end{diagram}

By $\BQ\text{HS}$ we denote the set of
rational homology 3--spheres $M$, with $|H_1(M,\Z)|=b$.
In particular, for any $M \in \BQ\text{HS}$,
we generalize Ohtsuki series as follows.
%For
% any rational homology 3--sphere $M$, with $|H_1(M,\Z)|=b$,
Let us fix a divisor $c<b$  of $b$,
%we construct power series of perturbative invariants
%dominating WRT invariants of $M$ at all roots of unity.
%Indeed,
then the Taylor expansion of $I_M$ at $e_c$ is a power series in $(q-e_c)$ with
coefficients in $\Z[1/b][e_c]$ which dominates the WRT invariants
at roots of unity whose order has the greatest common divisor $c$ with $b$.
If $c=b=p^k$, then  $I_{M,b}$ is a priory  determined by a product of power
series $\prod_{i\geq k} \BZ\left[\frac{1}{p}, e_{p^i}
\right]\left[\left[e_{{p^i}}-q\right]\right]$,
 however, conjecturally, it is enough to consider  the series in $q-e_b$.

The commutative diagram tells us that LMO invariant determines  $I_{M,1}$, or
the set of WRT invariants at roots of unity coprime with $b$.
Since there is no direct way to obtain power series in $(q-e_c)$ from the
LMO invariant, we conjecture the existence of the refined universal
finite type invariant, dominating our power series.
On the physical level of rigor, this means that the new refined
invariant should capture
more than just the contribution of flat connections
from the Chern--Simons theory.

The methods used in the unification
of the WRT invariants  also led to the full solution of
the integrality problem for quantum $SO(3)$ and $SU(2)$ WRT invariants.
In \cite{BL,BCL}, we showed that
$\tau_M^G(\xi)$ for any 3--manifold $M$ and any root of unity
is always an algebraic integer. Here $G=SO(3)$ or $SU(2)$.
The integrality of the  spin
and cohomological refinements is work in progress.

Assume $M$ is the Poincar\'e homology 3-sphere, obtained by surgery on a left-hand trefoil
with framing -1. Then
$$ I_M=\frac{1}{1-q} \sum^\infty_{k=0} q^k (1-q^{k+1})(1-q^{k+2})
\dots (1-q^{2k+1})\, .$$
We expect that  the categorification of
the WRT invariants will lead to  a
homology theory with  Euler characteristic given by $I_M$.

The paper is organized as follows. In Section \ref{defs} we recall
the definitions of Kirby moves,
WRT invariants and of the cyclotomic expansion for the colored Jones polynomial.
In Section 2 we state our main results and outline the proofs.
 Section 3 is
devoted to the discussion of the rings, where the unified
invariants take their values.
In addition, we construct generalized
Ohtsuki series.
The Laplace transform method, the Andrews identity, Frobenius maps
and the integrality of WRT invariants are explained in the last Section.
%We conjecture that $I_M$ belongs to a even smaller ring
%than shown in \cite{BBL}, such that the integrality of WRT invariants
%can be shown directly.
%\end{document}

%For rational homology 3--spheres the universal finite type
%invariant was constructed by
% Le, Murakami and Ohtsuki  \cite{LMO}. It
%In
% \cite{Le}, it is shown that the LMO invariant
% determines Ohtsuki series and, hence,
%$\{\tau_M(\xi)\,|\,(\ord(\xi),b)=1\}$.
%An interesting open question is whether
% the Le--Murakami--Ohtsuki invariant  determines
%$I_M$.

% Our
% work seems to show that the case  $(b,r) \neq 1$ is out of reach from  the
%original Ohtsuki series, and gives rise to new Ohtsuki's series.

%In this paper, the theory of perturbative 3--manifold invariants founds
%its incarnation.

\section{Quantum (WRT) invariants}
 \label{defs}
%Let us recall the definition of the WRT invariants.

\subsection{Notations and conventions}
We will consider $q^{1/4}$ as a free parameter. Let
\[
\{n\} = q^{n/2}-q^{-n/2},
 \quad  \{n\}!=
\prod_{i=1}^n \{i\} ,\quad  [n] =\frac{\{n\}}{\{1\}}.
% \quad
%\qbinom{n}{k} = \frac{\{n\}!}{\{k\}!\{n-k\}!}.
\]
We denote the set $\{1,2,3,\ldots\}$ by
$\N$.
We also use the following notation from $q$--calculus:
$$ (x;q)_n := \prod_{j=1}^n (1-x q^{j-1}).$$
Throughout this paper,
%unless stated otherwise,
  $\xi$ will be  a primitive root
of unity of {\em odd} order $r$ and $e_n:=\exp(2\pi I/n)$.

All 3--manifolds in this paper are supposed to be  closed and
oriented. Every link in $S^3$ is framed, oriented, and has
components ordered.
%We denote by $L$  a framed link in $S^3$  with
For a  link $L$,
let $(L_{ij})$ be its linking matrix, where for $i\neq j$,
$L_{ij}:=\lk(L_i,L_j)$ is the linking number of the $i$th and the $j$th
components,
and $L_{ii}=b_i$ is the
 framing of  $L_i$,
given by the linking number  of $L_i$ with its push off
along the first vector in the framing.

Surgery along the framed link $L$ consists of removing a tubular neighborhood of $L$
from $S^3$ and then gluing it back with the diffeomorphism given by framing
(i.e. the meridian of each removed solid torus, diffeomorhic to a neighborhood of $L_i$,
is identified with the
longitude of the complement twisted $b_i$ times along the meridian).
We denote by   $M$ the result of the surgery.
%L')$. We use the same notation $L'$ to denote the link in $S^3$ and
%the corresponding one in $M$.

%define framing!!

%%%%%%%%%%%%%%%%%%%%%%%%%%%%%%%%%%%%%%%%%%%%%%%%%%%%%%%%%%%%%%%%%%%%%%%%%%%%%%%%%%%

\subsection{The colored Jones polynomial}
\newcommand{\RR} {\mathbf R}

Suppose L is a framed oriented link with $m$ ordered components and $V_1,\dots,V_m$ are finite-dimensional modules
over a ribbon Hopf algebra. Then one can define the quantum invariant $J_L(V_1,\dots,V_m)$ through the machinery of quantum link invariant theory, see  \cite{Tu}.
The quantized enveloping algebra $U_h(sl_2)$ of $sl_2$ is a ribbon algebra, and
 for each positive integer $n$ there is
a unique $U_h(sl_2)$-module of dimension $n$. The quantum link invariant $J_L(n_1,\dots,n_m)$, where $n_j$ stands for the $n_j$-dimensional $U_h(sl_2)$-module, is usually called
the colored Jones polynomial, with the $n_j$'s being the colors. When all the colors are 2, $J_L(2,\dots,2)$ is the usual Jones polynomial, which can be defined using
skein relation.

One can construct the colored Jones polynomial without the quantum group theory  by first defining the Jones polynomial through the skein relation, and then
defining the colored Jones polynomial by using cablings, see e.g. \cite{KM, Lic}.

Let us recall here a few well--known formulas.
For the unknot $U$ with 0 framing one has
\begin{equation} J_U(n) = [n] \label{unknot}
\end{equation}
Moreover, $J_L$ is multiplicative with respect to the disjoint union.
 If $L_1$ is obtained from $L$
by increasing the framing of the $i$th component by 1, then
\begin{equation}\label{framing}
J_{L_1}(n_1,\dots,n_m) = q^{(n_i^2-1)/4} J_{L}(n_1,\dots,n_m).
\end{equation}
If all the colors $n_i$ are odd, then $J_{L}(n_1,\dots,n_m) \in \Z[q^{\pm 1}]$.

% \vskip2mm \noindent {\bf Exercise.} Prove \eqref{framing} by using one of the definitions of the colored Jones polynomial. \vskip2mm

%%%%%%%%%%%%%%%%%%%%%%%%%%%%%%%%%%%%%%%%%%%%%%%%%%%%%%%%%%%%%%%%%%%%%%%%%%%%%%%%%%%

\subsection{Kirby moves}
By Kirby theorem, any link invariant which does not change
under Kirby moves is an invariant of a 3--manifold given
by surgery on that link.
Let us first recall what the Kirby moves are.
\newcommand{\tL}{\tilde L}

\vskip2mm
{\noindent \bf K1--Move} (handle slide):
For some $i\neq j$, replace the $i$th component $L_i$
with  $L'_i$, a band connected sum of $L_i$ with a push
off of $L_j$ (defined by the framing), with
$b'_i=b_i+b_j+2\lk(L_i,L_j)$.

\vskip2mm
{\noindent \bf K2--Move} (blow up): Add (or delete) a split unknotted component with framing
$\pm 1$.
\vskip2mm

These two moves are equivalent to the one Fenn--Rourke move
defined as follows:
\vskip2mm
{\noindent \bf FR--Move} Locally the following two pictures are
interchangeable
$$\psdiag{10}{30}{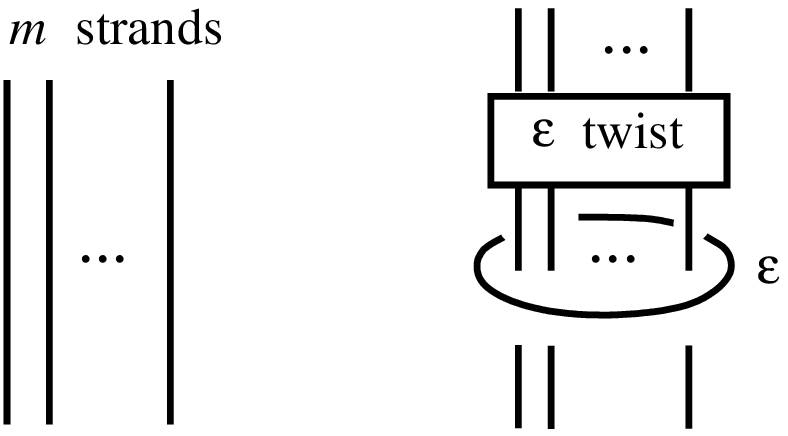}$$
where $\varepsilon \in \{1,-1\}$ and
the closed component has framing $\varepsilon$.
Note that K2--Move corresponds to the case when $m=0$.
\vskip2mm

The main idea of the construction of the WRT invariants
is to make the colored Jones polynomial invariant under Kirby moves
by  averaging over all colors.
To make this precise, we have to choose the quantum parameter $q$
to be a root of unity, otherwise the sum would be infinite.

\subsection{Evaluation and Gauss sums}
For each root of unity $\xi$ of odd order $r$,
we define the evaluation map
$\ev_\xi$ by replacing $q$ with $\xi$.
%Suppose $f(q; n_1,\dots, n_m)$
%is a function of variables $q^{\pm 1}$ and

% where $d$ is coprime with $r$. There
%exists an integer $d_*$, unique modulo $r$, such that
%$(\xi^{d_*})^{d}=\xi$. Then we define
%\[\ev_\xi f := f|_{q^{1/d}= \xi^{d_*}}.\]

Suppose
  $f(q;n_1,\dots,n_m)$ is a function
of variables $q^{\pm 1}$ and integers $n_1,\dots,n_m$.
In quantum topology, the following sum plays an important role
\[
{\sum_{n_i}}^\xi f := \sum_{\substack{0< n_i< 2r\\ n_i \text{ odd}}}
\ev_\xi  f(q; n_1,\dots, n_m)
\]
where in the sum  all the $n_i$ run over the set of {\em odd} numbers
between $0$ and $2r$.

In particular, the following  variation  of the Gauss sum
\[
\gamma_b(\xi):= {\sum_{n}}^\xi q^{b\frac{n^2-1}{4}}
\]
is well--defined, since for odd $n$, $4\mid n^2-1$.
It is known that, for odd $r$, $|\gamma_b(\xi)|$ is never 0.

\subsection{Definition of the WRT invariants}
%Suppose the components of $L'$ are colored by fixed integers
%$j_1,\dots,j_l$.
 Let
\[
%F_L(\xi):= {\sum_{n_i}}^\xi
%J_L(n_1,\dots,n_m)\prod_{i=1}^m [n_i]\, ,
 F_{L}(\xi):= {\sum_{n_i}}^\xi\; \prod_{i=1}^m [n_i]\;
  J_{L}(n_1,\dots,n_m).\]

\begin{Thm} \label{K1M}{\rm [Reshetikhin--Turaev]}
$F_L(\xi)$ is invariant under K1--Move.
\end{Thm}

% \vskip2mm \noindent {\bf Exercise.} Prove this theorem by using  FR--Move with $m>0$, fusion rules, \eqref{framing} and the fact that the Jones polynomial of the $0$--framed $(i,j)$--colored Hopf link is $[ij]$. (Hint: use Lemma \ref{33} below.) \vskip2mm

An important special case is when $L=U^b$, the unknot with framing
$b \neq 0$. In that case $F_{U^{b}}(\xi)$ can be calculated
using the Gauss sum and is nonzero.

 Let $\sigma_+ $ (respectively $\sigma_-$) be the number of
positive (negative) eigenvalues of the linking matrix of $L$.
 Then we define
%the quantum
%$SO(3)$ invariant of the pair $(M, L')$ is defined by (see e.g. \cite{KM,Tu})
\begin{equation}
\tau_{M}(\xi) =
\frac{F_{L}(\xi)}{(F_{U^{+1}}(\xi))^{\sigma_+}\,
(F_{U^{-1}}(\xi))^{\sigma_-} }\, .
\label{def_qi}
\end{equation}
%where $U^\pm$ is the unknot with framing $\pm 1$.
It is easy to see that  $\tau_M(\xi)$ is invariant under
K2--Move, and hence, by Theorem \ref{K1M}, it is a topological
invariant of $M$ called the $SO(3)$ WRT invariant.
 Moreover,
 $\tau_{M}(\xi)$ is multiplicative with
respect to the connected sum.

\vskip2mm
\noindent
{\bf Remark.}
If we drop the condition that the colors in the summation are
odd, and sum over all (odd and even) colors, the result will be
the $SU(2)$ WRT invariant  $\tau^{SU(2)}_M (\xi)$.
In this case, the order of $\xi$ could also be even.
\vskip2mm

 The $SO(3)$  WRT invariant  extends naturally to the invariant
of the pair $(M,L')$, where
the manifold $M$ contains a link $L'$ inside. In this case we have
to replace the surgery link $L$ of $M$ by $L\cup L'$ in all definitions,
fix colors on $L'$ and sum over all colorings of
$L$ only. We omit here the precise definition and refer to \cite{BBL}
for more details.

For example, the $SO(3)$ invariant of the lens space $L(b,1)$,
obtained by surgery along $U^b$,  is
\begin{equation} \tau_{L(b,1)} (\xi)= \frac{ F_{U^b}(\xi)}{F_{U^{\sn(b)}}(\xi) },
\label{2005}
\end{equation}
were $\sn(b)$ is the sign of the integer $b$.

Suppose
$M$ is a rational homology
3--sphere.
%$|H_1(M,\Z)|:={\rm card}\,H_1 (M,\Z) < \infty$.
Then there is a unique decomposition
 $ H_1(M,\Z)=\bigoplus_{i}
\Z/{b_{i}\Z}$, where each $b_i$ is a prime power.
We renormalize the $SO(3)$  WRT invariant  as follows:
\begin{equation}
\tau'_{M}(\xi)=\frac{\tau_{M}(\xi)}
{\prod\limits_{i}\;\,
\tau_{L(b_{i},1)}(\xi)}\; ,
%\;\;\;\;\;\;\;\;\;\;
%\tau'^{SU(2)}_{M}(\xi)=\frac{\tau^{SU(2)}_M(\xi)}
%{\prod\limits_{i}\;\,
%\tau^{SU(2)}_{L(p^{k_i}_{i},1)}(\xi)}\; ,
\label{0910}
\end{equation}
where $L(b,a)$ denotes the $(b,a)$ lens space. Note that
$\tau_{L(b,1)}(\xi)$ is always nonzero.

Let us focus on the special case when the linking matrix of
$L$ is diagonal.
%with $b_1, b_2, \dots, b_m$ on the diagonal.
Assume each $L_{ii}=b_i$ is a power of a prime or 1, up to sign.
Then $H_1(M,\Z) = \oplus_{i=1}^m \Z/|b_i|$, and
$$\sigma_+ = {\rm card}\, \{ i\mid b_i >0\}, \quad \sigma_- =
{\rm card}\, \{ i \mid b_i < 0\}.$$
Thus from the definitions \eqref{def_qi}, \eqref{2005} and \eqref{0910} we have
\begin{equation}
 \tau'_{M}(\xi) = \left( \prod_{i=1}^m  \tau'_{L(b_i,1)}(\xi)
 \right)\,
\frac{F_{L}(\xi)}
{\prod_{i=1}^m F_{U^{b_i}}(\xi) }
\,  ,
\label{0077}
\end{equation}
with
$$\tau'_{L(b_i,1)}(\xi)= \frac{\tau_{L(b_i,1)}(\xi)}{\tau_{L(|b_i|,1)}(\xi)}\, .$$

The collection of $SO(3)$ WRT invariants is difficult to study, since
their definition heavily depends on the order of $\xi$.
The following cyclotomic expansion will play an important role
for the unification of the WRT invariants.

\subsection{Habiro's cyclotomic expansion of the colored Jones
polynomial}
%Recall that  $L$ and $L'$ have $m$ and $l$ components, respectively. Let us color $L'$ by fixed $\bj=(j_1,\dots,j_l)$ and vary the colors
%$\bn=(n_1,\dots,n_m)$ of $L$.
For non--negative integers $n,k$ we define
$$ A(n,k) := \frac{\prod^{k}_{i=0}
\left(q^{n}+q^{-n}-q^i -q^{-i}\right)}{(1-q) \, (q^{k+1};q)_{k+1}}.$$
For $\bk=(k_1,\dots,k_m)$ let
$$ A(\bn,\bk):= \prod_{j=1}^m \; A(n_j,k_j).$$
Note that $A(\bn,\bk)=0$ if $k_j \ge n_j$ for some index $j$. Also
$ A(n,0)= q^{-1} J_U(n)^2.$

The colored Jones polynomial $J_{L}
(\n)$  can be repackaged into the invariant $C_{L}
(\bk)$ as stated in the following theorem.
\begin{Thm} {\rm [Habiro]}\label{GeneralizedHabiro}
 Suppose $L$ is a link in $S^3$ having zero linking matrix.
%Assume the components of $L'$ have fixed {\em odd}
%colors $\bj = (j_1,\dotsm j_l)$.
Then  there are invariants
\begin{equation}\label{Jones2}
C_{L}(\bk) \in \frac{(q^{k+1};q)_{k+1}}{(1-q)}
\,\,\BZ[q^{\pm 1}] ,\quad \text{where  $k=\max\{k_1,\dots, k_m\}$}
\end{equation}
such that for every $\bn =(n_1,\dots, n_m)$
\begin{equation}\label{Jones}
 J_{L}
(\n)  \,  \prod^m_{i=1}\; [n_i] = \sum_{0\le k_i \le n_i-1}
C_{L}(\bk)\;
  A(\bn, \bk).
\end{equation}
\end{Thm}

% When $L'=\emptyset$,
%This is
% Theorem 8.2 in \cite{Ha}.
 Note that the existence of $C_{L}(\bk)$ as
rational functions in $q$ satisfying \eqref{Jones} is
 easy to establish. They correspond to the Jones polynomial
colored by different elements of the Grothendieck ring
of $U_q(\frak {sl}_2)$, i.e. by linear combinations of representations.
The difficulty here is to show the integrality
of \eqref{Jones2}.

 Since $A(\bn, \bk) =0$ unless $ \bk < \bn$, in the sum on the right
hand side of \eqref{Jones} one can assume that $\bk$ runs over the set
 of all $m$--tuples $\bk$ with non--negative integer components. We will use
this fact later.

% Suppose
%$L$ is an algebraically split link with 0--framing on each component.
%Then we have
%\begin{equation}\label{evJones}
%\ev_\xi (J_{L\sqcup L'}(n_1,\dots,n_m, j_1,\dots,j_l))  = \ev_\xi \left(
%\sum_{k_1,\dots,k_m=0}^{(r-3)/2} J_L(P'_{k_1}, \dots , P'_{k_m},
%V_{j_1},\dots, V_{j_l})
%\prod_{i=1}^m \qbinom{n_i+k_i}{2k_i+1} \{k_i\}!\right).
%\end{equation}

%%%%%%%%%%%%%%%%%%%%%%%%%%%%%%%%%%%%%%%%%%%%%%%%%%%%%%%%%%%%%%%%%%%%
%%%%%%%%%%%%%%%%%%%%%%%%%%%%%%%%%%%%%%%%%%%%%%%%%%%%%%%%%%%%%%%%%%%%%%%%%%%%%%%%%%%%%%%%%%%%%%%%%%%%%%%%%%%%%%%%%%%%%%%%%

\section{Main Results}
\label{strategy}
Let us state our main results announced in the introduction.

%Suppose
%$M$ is a rational homology
%3--sphere  with a unique decomposition
% $ H_1(M,\Z)=\bigoplus_{i}
%\Z/{b_{i}\Z}$, where each $b_i$ is a prime power and $b=\prod_i b_i$.

%\begin{equation}
%\tau'_{M, L}(\xi)=\frac{\tau_{M,L}(\xi)}
%{\prod\limits_{i}\;\,
%\tau_{L(b_{i},1)}(\xi)}\; ,
%\;\;\;\;\;\;\;\;\;\;
%\tau'^{SU(2)}_{M}(\xi)=\frac{\tau^{SU(2)}_M(\xi)}
%{\prod\limits_{i}\;\,
%\tau^{SU(2)}_{L(p^{k_i}_{i},1)}(\xi)}\; ,
%\label{0910}
%\end{equation}
%where $L(b,a)$ denotes the $(b,a)$ lens space. We will see that
%$\tau_{L(b,1)}(\xi)$ is always nonzero.

 For any positive integer $b$, we  define the
cyclotomic completion ring $\R_b$ to be
\be
\label{ab} \R_b:=\lim_{\overleftarrow{\hspace{2mm} k\hspace{2mm}}}
\frac{\Z[1/b][q]}{
\left((q;q^2)_k\right)}, \qquad \text{where} \quad
(q;q^2)_k = (1-q)(1-q^3) \dots (1-q^{2k-1}).
\ee
For any    $f(q)\in \R_b$ and a root of unity $\xi$ of {\em odd} order,
the evaluation $\ev_\xi (f(q)):= f(\xi)$ is well--defined.
Similarly, we put
$$\calS_b :=\lim_{\overleftarrow{\hspace{2mm}k\hspace{2mm}}}
\frac{\Z[1/b][q]}
{((q;q)_k)}\; .$$
Here the evaluation at any root of unity is well--defined. For odd $b$,
there is  a natural embedding
$\calS_b\hookrightarrow \R_b$.

Let us denote by $\M_b$ the set of rational homology 3--spheres such that
$|H_1(M,\Z)|$ divides $b^n$ for some $n$.
Our main result is the following.
\begin{Thm}\label{main} {\rm [Beliakova--B\"uhler--Le]}
Suppose the components of a framed oriented link $L \subset M$ have odd colors, and $M\in \cM_b$. Then
there exists an invariant $I_{M,L} \in \R_b$,
such that for any root of unity $\xi$ of odd order
$$\ev_\xi(I_{M,L})=\tau'_{M,L}(\xi)\, .$$
In addition, if $b$ is odd, then
$I_{M,L}\in \calS_b$.
\end{Thm}

If $b=1$ and $L$ is the empty link, $I_{M}$
coincides with Habiro's unified invariant $J_M$ and
$\calS_1=\Habiro$.

One may ask what is the evaluation of $I_M$ at an even
root of unity in the case when $b$ is odd.
In her PhD thesis \cite{B}, B\"uhler shows that it
coincides with $\tau^{\prime SU(2)}_M(\xi)$.
Hence, for odd $b$, $I_M$ dominates both $SO(3)$ and $SU(2)$
WRT invariants. An analogous result for $b$ even is
work in progress.

Compared to  Habiro's case,
the proof of Theorem \ref{main} uses the following
new techniques:
1)  the Laplace transform method;
2) the difficult number theoretical identity of Andrews generalizing those of
Roger--Ramanujan;
3) the Frobenius type isomorphism providing the existence of the
$b$--th root of $q$ in $\R_b$.
%We also construct a Frobenius type isomorphism to get rid of the
%formal fractional power of $q$ that appeared in \cite{Le}, \cite{BL}.
In addition,  we had to generalize the deep integrality result of Habiro
(Theorem \ref{GeneralizedHabiro}),
  to a union of an algebraically split link
with any odd colored one.

The rings $\R_b$ and $\calS_b$ have properties
 similar  to those of the Habiro ring.
An element $f(q) \in \R_b$ is totally determined by the values at many
infinite sets of roots of unity (see Section \ref{cyc}),
one special case is the following.

\begin{proposition}\label{main-cor} {\rm [Beliakova--B\"uhler--Le]}
Let $p$ be an odd prime
not dividing $b$ and $T$ the set of all integers of the form
$p^k b'$ with $k\in \N$ and $b'$ any odd divisor of
$b^n$ for some $n$. Any element $f(q) \in \R_b$, and hence also
$\{\tau_M(\xi)\}$, is totally determined by the values at roots of
unity with orders in $T$.
\end{proposition}

The general properties of the ring $\R_b$
allow to introduce generalized
 Ohtsuki series as the Taylor expansions
of $I_M$ at roots of unity. In addition, we show
 that these Taylor expansions satisfy
congruence relations similar to the original definition of the Ohtsuki
series (see Section \ref{Oht}).

\subsection{Strategy of the proof}
Let us outline the proof of Theorem \ref{main} and state the main
technical results that will be explained later.

%As before, $ L\sqcup L'$ is a framed link in $S^3$ with
%disjoint sublinks $L$ and $L'$, with $m$ and $l$ components, respectively.
%Assume that $L'$ is colored by fixed $\bj=(j_1,\dots, j_l)$, with $j_i$'s odd.
%Surgery along the framed link $L$ transforms $(S^3,L')$ into  $(M,
%L')$.

We restrict to the case $L=\emptyset$ for simplicity.
We would like to define $I_{M}\in \cR_{b}$, such that
\begin{equation}
\tau'_{M}(\xi)\;=\; \ev_{\xi}\left(I_{M}\right)
\label{0080}
\end{equation}
for any root of unity $\xi$ of odd order.
This unified invariant is multiplicative with respect
to the connected sum.

The following observation is important.
By Proposition \ref{main-cor}, there is {\em at most one} element $f(q)\in \R_b$  such that for every root $\xi$ of odd order one has
$$ \tau'_{M} (\xi) = \ev_\xi\left( f(q)\right).$$
That is, if  we can find such an element,
 it is unique, and we put
$I_{M} := f(q)$.

\subsection{Laplace transform}
The following is the main technical result of \cite{BBL}. A proof will be
explained in the next Section.
\begin{Thm} {\rm [Beliakova--B\"uhler--Le]}
 Suppose $b=\pm 1$ or $b= \pm p^l$ where $p$  is a  prime and $l$ is
positive. For  any non--negative integer $k$,
there exists an element $Q_{b,k} \in \R_b $ such that
for every root $\xi$ of odd order $r$ one has
\[
\frac{{\sum\limits_n}^\xi \, q^{b\frac{n^2-1}{4}} A(n,k) }{F_{U^b}(\xi)}
= \ev_\xi (Q_{b,k}).
\]
\label{0078}
In addition, if $b$ is odd, $Q_{b,k} \in \calS_b $.
\end{Thm}

% For $\bb:= (b_1,\dots, b_m)$ and $\bk:=(k_1,\dots,k_j)$ let $$ Q(\bb,\bk) := \prod_{i=1}^m Q(b_i,k_i).$$

\subsection{Definition of the unified invariant: diagonal case} \label{2501}
Suppose that the linking number between any two components of $L$
is 0, and  the framing on components of $L$ are $b_i=\pm p_i^{k_i}$ for
$i=1,\dots, m$, where each $p_i$ is prime or 1.
Let us denote the link $L$ with all framings switched to zero by $L_0$.
%Further suppose that $c$ is a fixed odd divisor of $b=\prod_{i} b_{i}$,
%$c_{i}=(c,b_{i})$ and $t_{i}$ is a primitive $b_i/c_i$th root of $q^{c_i}$.

Using \eqref{Jones}, taking into account the framings $b_i$'s, we have
\[
J_{L}(\bn)\prod_{i=1}^m [n_i] = \sum_{\bk\ge 0} C_{L_0 }
(\bk) \, \prod_{i=1}^m q^{b_i \frac{n_i^2-1}{4}} A(n_i,k_i).
\]
By the definition of $F_{L}$, we have
\[ F_{L}(\xi)= \sum_{\bk \ge 0}
\ev_\xi(C_{L_0 }(\bk)) \, \prod_{i=1}^m
{\sum_{n_i}}^\xi \,  q^{b_i \frac{n_i^2-1}{4}} A(n_i,k_i).
\]
From \eqref{0077} and  Theorem \ref{0078}, we get
\[
\tau'_{M}(\xi) = \ev_\xi \left \{  \prod_{i=1}^m  I_{L(b_i,1)} \,
 \sum_{\bk} C_{L_0 }(\bk) \, \prod_{i=1}^m  Q_{b_i,k_i} \right \},
\]
where the existence of
 the unified invariant of the lens space $I_{L(b_i,1)}\in \R_b$,
with
 $\ev_\xi(I_{L(b_i,1)})=\tau'_{L(b_i,1)}(\xi)$ can be shown by a
direct computation (we refer to \cite{BBL} for more details).

Thus if we define
\[
I_{M}:= \prod_{i=1}^m  I_{L(b_i,1)} \, \sum_{\bk}
C_{L_0 }(\bk) \, \prod_{i=1}^m  Q_{b_i,k_i}\, ,
\]
then \eqref{0080} is satisfied. By Theorem \ref{GeneralizedHabiro},
$C_{L_0 }(\bk)$ is divisible by $(q^{k+1};q)_{k+1}/(1-q)$, which is
divisible by $(q;q)_k$, where $k = \max k_i$. It follows that
$I_{M} \in \R_b$. In addition, if $b$ is odd,
then $I_{M} \in \calS_b$.

\subsection{Diagonalization using lens spaces} The general
case reduces to the diagonal case by
the well--known trick of diagonalization using lens spaces. We say that
$M$ is {\em diagonal} if it can be obtained from $S^3$ by surgery along a framed link
$L$  with diagonal linking matrix, where  the diagonal entries are of the
form $\pm p^k$ with $p=0,1$ or a prime.
The following lemma was proved in \cite[Proposition 3.2 (a)]{Le}.
% we will give a proof in Section
%\ref{diag}. For various
% versions see \cite{Le,Oh}.

\begin{lemma} For every rational homology sphere $M$,
 there are lens spaces $L(b_i,a_i)$ such that the connected
sum of $M$ and these
lens spaces is diagonal. Moreover, each $b_i$ is a prime power
%(i.e. $b_i=p^k$ for some prime $p$)
 divisor of $|H_1(M,\Z)|$.
\label{diagonalization}
\end{lemma}

To define  the unified invariant for a general rational homology sphere $M$,
one first adds to $M$ lens spaces to get a diagonal $M'$, for which
the unified invariant $I_{M'}$ had been defined in Subsection \ref{2501}. Then
  $I_M$ is the quotient of $I_{M'}$ by the unified invariants of the lens
 spaces. But unlike the simpler
 case
of \cite{Le}, the unified invariant of lens spaces are {\em not} invertible
 in general.
To overcome this difficulty we insert knots in lens spaces and split the
 unified invariant into different components.
This is also the reason why we need to generalize Habiro's integrality
result to algebraically split links together with odd colored components.

%This
% is  explained in \cite{BBL} in details.

%%%%%%%%%%%%%%%%%%%%%%%%%%%%%%%%%%%%%%%%%%%%%%%%%%%%%%%%%%%%%%%%%%%%%%%%%%%%%%%%%%%%%%%%%%%%%%%%%%%%%%%%%%%%%%%%%%%%%%%%%%%%%%%%%%%%%%%%%%%%%%%%%%%%%%%%%%%%%%%%%%%%%%%%%%%%%%%%%%%%%%%%%%%%%%%

\section{Cyclotomic completions of polynomial rings}\label{cyc}
Since unified invariants belong to
 cyclotomic completions of polynomial rings, we outline
their construction and the main properties. For simplicity, only the case
$b$ is a power of a prime is considered, the general case is treated in
\cite{BBL}.

\subsection{On cyclotomic polynomials}
Recall that $e_n := \exp(2\pi I/n)$
and denote by
$\Phi_n(q)$ the cyclotomic polynomial
\[
\Phi_n(q) = \prod_{\substack{(j,n)=1\\0<j<n}} (q - e_n^j).
\]
For example, $\Phi_1(q)=q-1$ and $\Phi_2(q)=q+1$.
The degree of $\Phi_n(q)\in \Z[q]$ is given by the Euler function $\varphi(n)$.
Suppose $p$ is a prime and $n$ an integer. Then (see e.g. \cite{Na})
\begin{equation} \Phi_n(q^p)= \begin{cases}    \Phi_{np}(q)  & \text{ if } p \mid n \\
\Phi_{np}(q) \Phi_n(q)  & \text{ if } p \nmid n.
\end{cases}
\end{equation}
It follows that  $\Phi_n(q^p)$ is always divisible by $\Phi_{np}(q)$.

The ideal of $\Z[q]$ generated by $\Phi_n(q)$ and $\Phi_m(q)$ is well--known,
see e.g. \cite[Lemma 5.4]{Le}:

\begin{lemma}$\text{ }$
\begin{itemize}
\item[(a)] If $\frac{m}{n} \neq p^e$ for any  prime $p$ and any integer $e\neq 0$, then
$(\Phi_n)+ (\Phi_m)=(1)$ in $\Z[q]$.

\item[(b)] If $\frac{m}{n} = p^e$ for a prime $p$ and some integer $e \neq 0$, then $(\Phi_n)+ (\Phi_m)=(p)$ in $\Z[q]$.
\label{0911}
\end{itemize}
\end{lemma}

%Note that in a commutative ring $R$,  $(x) + (y) =(1)$ if and only if
%$x$ is invertible in $R/(y)$. Also $(x) + (y) =(1)$ implies $(x^k) +
%(y^l) =(1)$ for any integers $k,l \ge 1$.

%\begin{lemma}Suppose  $r$  is coprime with $p$, and $x,y \in\Q[e_r]$
% such that$x^k=y^k$ for some $k>1$ which is a power of $p$. Then $x=y$.
%\label{1102}
%\end{lemma}
%
%\begin{proof} Then $x/y$ is a root of 1 of order a power of $p$.
%If $x/y\neq 1$, then this means $\Q[e_r]$ contains a
%primitive $p$th root of 1, or $e_p \in \Q[e_r]$, which is %impossible.
%\end{proof}

%%%%%%%%%%%%%%%%%%%%%%%%%%%%%%%%%%%%%%%%%%%%%%%%%%%%%%%%%%%%%%%%%%%%%%%%%%%%%%%%%%%%%%%%%%%%%%%

\subsection{Habiro's results}  Let us summarize some of Habiro's
results on cyclotomic completions of polynomial rings \cite{Ha1}.
Let  $R$ be a commutative integral domain of characteristic zero
and $R[q]$  the polynomial
ring over $R$.
%For each $n\in \N$, let
%\[
%\Phi_n(q): =\prod_{(i,n)=1} (q-e_n^i)
%\]
%denotes the $n$th cyclotomic polynomial, where $e_n$ is a primitive $n$th
%root of unity.
For any
 $S\subset \N$, Habiro defined  the $S$--cyclotomic completion
 ring $R[q]^S$ as follows:
\be\label{rs} R[q]^S:=\lim_{\overleftarrow{f(q)\in
\Phi^*_S}} \;\;\frac{R[q]}{(f(q))} \ee
 where $\Phi^*_S$ denotes the multiplicative
set in $\Z[q]$ generated by $\Phi_S=\{\Phi_n(q)\mid n\in S\}$
  and directed with respect to the
divisibility relation.

For example, since the sequence $(q;q)_n$, $n\in \N$,
 is cofinal to $\Phi^*_\N$, we have
\be\label{cofinal}
\Habiro\simeq\Z[q]^\N.
\ee

 Note that if $S$ is finite, then
$R[q]^S$ is identified with the $(\prod \Phi_S)$--adic completion of $R[q]$.
In particular,
\[
R[q]^{\{1\}}\simeq R[[q-1]], \quad
R[q]^{\{2\}}\simeq R[[q+1]].
\]
Suppose $S' \subset S$, then $\Phi^*_{S'}\subset \Phi^*_S$, hence
there is a natural map
\[
\rho^R_{S, S'}: R[q]^S \to R[q]^{S'}.
\]

Recall important results concerning  $R[q]^S$ from \cite{Ha1}. Two
positive integers $n, n'$ are called {\em adjacent} if
$n'/n=p^e$ with a nonzero
$e\in \Z$ and  a  prime $p$, such that
 the ring $R$ is $p$--adically separated, i.e. $\bigcap_{n=1}^\infty (p^n) =0$ in $R$.
A set of positive integers is
{\em $R$--connected} if for any two distinct elements $n,n'$ there is a
sequence $n=n_1, \,n_2, \dots,\, n_{k-1},\, n_k= n'$ in the set, such that
any two consecutive numbers of this sequence are adjacent.
Theorem 4.1 of \cite{Ha1}
says that if $S$ is $R$--connected,
 then for any subset $S'\subset S$
the natural map $ \rho^R_{S,S'}: R[q]^S \hookrightarrow R[q]^{S'}$ is an
embedding.

\vskip2mm
\noindent
{\bf Example}: Assume $R=\Z$, $S=\N$ and $S'=\{1\}$, then
we have that the map
$\Habiro \to \Z[[q-1]]$ is an embedding. This implies that,
for any integral homology 3--sphere $M$,
$J_M$ is determined by the Ohtsuki series and this series has
integral coefficients.
\vskip2mm

If $\zeta$ is a root of unity of order in $S$, then for every $f(q)\in R[q]^S$
the evaluation $\ev_\zeta(f(q))\in R[\zeta]$ can be defined by sending
$q\to\zeta$.
For a set $\Xi$ of roots of unity whose orders form a subset
$\cT\subset S$, one defines the evaluation
\[
\ev_\Xi: R[q]^S \to \prod_{\zeta \in \Xi} R[\zeta].
\]
Theorem 6.1 of \cite{Ha1} shows that if
$R\subset \Q$, $S$ is $R$--connected, and there
exists $n\in S$ that is adjacent to infinitely many elements in
$\cT$, then $\ev_\Xi$ is injective.

\vskip2mm
\noindent
{\bf Example}: Consider again the case when $R=\Z$, $S=\N$
and put $\cT=\{3^n| n\in \N\}$, then $3\in S$ is adjacent to infinitely many
elements of $\cT$ and hence, for any integral homology 3--sphere $M$,
the whole set of its WRT invariants is determined by the evaluations
of $J_M$ at roots of unity of order in $\cT$.
\vskip2mm

%%%%%%%%%%%%%%%%%%%%%%%%%%%%%%%%%%%%%%%
%%%%%%%%%%%%%%%%%%%%%%%%%%%%%%%%%%%%%%%%

\subsection{Taylor expansion} Fix a natural number $n$, then we have
$$ R[q]^{\{n\}} =
\lim_{\overleftarrow{\hspace{2mm}k\hspace{2mm}}}
\frac{R[q]}{(\Phi^k_n(q))}\; .$$
Suppose $ \Z \subset R \subset \Q$, then the natural algebra homomorphism

$$ h: \frac{R[q]}{(\Phi^k_n(q))} \to \frac{R[e_n][q]}{((q-e_n)^k)}$$
can be proved to be
 injective. Taking the inverse limit, we see that there is a natural injective algebra homomorphism

$$ h : R[q]^{\{n\}} \to R[e_n][[q-e_n]].$$

Suppose $n \in S$. Combining $h$ and $\rho_{S, \{n\}}: R[q]^S \to R[q]^{\{n\}}$, we get an algebra map

\newcommand{\TT}{{\mathfrak t}}

$$ \TT_n: R[q]^S \to R[e_n][[q-e_n]].$$
If $f\in R[q]^S$, then $\TT_n(f)$ is called the Taylor expansion of $f$ at $e_n$.

\subsection{Splitting of $\calS_p$ and evaluation}
For every integer $a$, we put $\N_a := \{ n \in \N \mid (a,n)=1\}$.

Suppose $p$ is a prime. Analogously to \eqref{cofinal},
we have
$$\calS_p \simeq \Z[1/p][q]^\N\,.$$

Observe that $\N$ is not $\Z[1/p]$--connected. In fact one has
 $\N =\amalg_{j=0}^\infty \; p^j \N_p$, where each $p^j\N_p$ is
$\Z[1/p]$--connected.
Let us define
\[
\calS_{p,j}:= \Z[1/p][q]^{p^j \N_p}.
\]
Note that for every $f \in \calS_p$, the evaluation $\ev_\xi(f)$ can be
 defined for every root $\xi$ of unity. For $f\in \calS_{p,j}$, the evaluation
$\ev_\xi(f)$ can be defined  when $\xi$ is a root of unity of order
in $p^j\N_p$. In \cite{BBL} we proved the following.

\begin{proposition} For every prime $p$ one has
\begin{equation} \calS_p\simeq \prod_{j=0}^\infty \calS_{p,j}.
\label{0912}\end{equation}

\end{proposition}

Let $\pi_j: \calS_p \to \calS_{p,j}$ denote the projection onto the
$j$th component in the above decomposition.
Suppose $\xi$ is a root of unity of order $r= p^j r'$,
with $(r',p)=1$. Then
for any $x\in \calS_p$, one has
$$ \ev_\xi(x) = \ev_\xi(\pi_j(x)).$$
If $i\neq j$ then $\ev_\xi(\pi_i(x))=0$.

%({\bf Splitting of $\R_p$, of $\calS_b$}??)

%%%%%%%%%%%%%%%%%%%%%%%%%%%%%%%%%%%%%%%%%%%%%%%%%%%%%%%%%%%%%%%%%%%%%%%%%%%%%%%%%%%%%%%%%%%%%%%%%%%%%%%%%%%%%%%%%%%%%%%%%%%%%%%%%%%%%%%%%%%%%%%%%%%%%%%%%%%%%%%%%%%%%%%%%%%%%%%%%%%%%%%%%

\subsection{On the  Ohtsuki series at roots of unity}\label{Oht}
 Suppose $M$ is a rational
homology 3--sphere with $|H_1(M,\Z)|=b$.
By Theorem \ref{main}, for any root of unity $\xi$ of order $pr$
$$\tau'_M(\xi)\in \Z[1/b][e_{pr}] \simeq \frac{\Z[1/b,e_r][x]}{(f_p(x+e_r))}\; .$$ where
\[
f_p(t):= \frac{t^p - e_r^p}{t-e_r}.
\]
Hence we can write
\be\label{Ohtsuki}
\tau'_M(e_{r}e_{p})= \sum_{n=0}^{p-2} a_{p,n} x^n
\ee
where  $a_{p,n}\in \Z[1/b,e_r]$.
The following proposition proven in \cite{BBL}
shows that the coefficients $a_{p,n}$
stabilize as $p\to \infty$.
\begin{proposition}{\rm [Beliakova--B\"uhler--Le]}
\label{main-cor1} Suppose $M$ is a
rational homology 3--sphere with $|H_1(M,\Z)|=b$, and $r$ is an odd
positive integer.
For every non--negative integer $n$,
there exists a unique invariant $a_n= a_n(M) \in \Z[1/b,e_r]$ such that
for every prime $p > \max (b,r)$,
 we have
\begin{equation}
a_n\equiv a_{p,n} \pmod p\;\;\; \text{in $ \Z[1/b,e_r]$ for} \;\;\;
0\le n \le  p-2.
\label{5501}
\end{equation}
Moreover, the formal series $\sum_{n} a_n (q-e_r)^n$ is equal to the Taylor
 expansion of the unified invariant $I_M$ at $e_r$.
\end{proposition}

%%%%%%%%%%%%%%%%%%%%%%%%%%%%%%%%%%%%%%%%%%%%%%%%%%%%%%%%%%%%%%%%%%%%%%%%%%%%%%%%%%%%%%%%%%%%%%%%%%%%%%%%%%%%%%%%%%%%%%%%%%%%%%%%%%%%%%%%%%%%%%%%%%%%%%%%%%%%%%%%%%%%%%%%%%%%%%%%%%%%%%%%%%%%%%%

\section{Laplace transform, Andrews identity
 and Frobenius maps} \label{map-frob}
The aim of this section is to define the Laplace transform
 and to study its image.

\subsection{Laplace transform}
To define the unified invariant we have to compute
$${\sum_n}^\xi q^{b\frac{n^2-1}{4}} A(n,k)$$
where  terms depending on $n$ in $A(n,k)$ look as follows
$$\prod^k_{i=0} (q^n+q^{-n}-q^i-q^{-i})\,=(-1)^{k+1}\;(q^n;q)_{k+1}\;
(q^{-n};q)_{k+1}\, .$$
Formally, the last expression can be considered as a Laurent polynomial
in $q$ and $q^{n}$. Hence we only have to compute
$${\sum_n}^\xi q^{b\frac{n^2-1}{4}} q^{na}\; ,$$
which can be easily done by the square completing argument.
Let us state the result.

Suppose $r$ is an odd number,  $b$ is positive integer and
$$ c:= (r,b), \quad b_1:= b/c, \quad r_1:=r/c.$$
\newcommand{\Tor}{{\rm Tor}}
\begin{lemma}{\rm [Beliakova--Le]}
\label{33}
One has
\begin{equation}
 \gamma_b(\xi) = c \,\gamma_{b_1}(\xi^c). \label{02}
 \end{equation}
\begin{equation}{\sum_{n}}^\xi q^{b \frac{n^2-1}{4}}\, q^{an} = \begin{cases}
0
&\text{if $c\, \nmid \, a$;}\\
 (\xi^c)^{-a_1^2 b_1^*}\; \gamma_b(\xi)  & \text{if $a=ca_1$},
 \end{cases}
 \end{equation}
where $b_1^\ast$ is an integer satisfying $b_1 b_1^* \equiv 1 \pmod{r_1}$.
\end{lemma}

This computation inspired us to introduce the following operator,
called the Laplace transform. Remember that $\int e^{ax^k} f(x)\, dx$
is called the Laplace transform
of $f$ of order $k$.

Let $\cL_{b,c;n}: \BZ[q^{\pm n},q^{\pm 1}] \to \BZ[q^{\pm  c/b}]$ be
the $\BZ[q^{\pm 1}]$--linear operator, called the (discrete)
Laplace transform (of the second order),
defined by
\be\label{four}
\cL_{b,c;n}(q^{na}) := \begin{cases}
0
&\text{if $c\, \nmid\,  a$;}\\
q^{-a^2/b}   & \text{if $a=ca_1$},
 \end{cases}
\ee
\begin{lemma} {\rm [Beliakova--Le]}
Suppose  $f \in \BZ[q^{\pm n},q^{\pm 1}]$. Then
$$ {\sum_{n}}^\xi q^{b\frac{n^2-1}{4}}  f = \gamma_b(\xi) \,
\ev_\xi(\cL_{b,c;n}(f)).
$$
\label{1000}
\end{lemma}

The point is that $\cL_{b,c;n}(f)$, unlike the left hand side
${\sum_{n}^\xi}  q^{b\frac{n^2-1}{4}} f$, does not depend on $\xi$.
To prove Theorem \ref{0078} we need  to show that
$\cL_{b,c;n} ((q^{n};q)_{k+1}\; (q^{-n};q)_{k+1})$ is divisible by
$(q^{k+1};q)_{k+1}$.
 For  this we use the remarkable identity
discovered by Andrews.
%and will help us to define a ``universal invariant".

\subsection{Andrews identity}
To warm up
we start with the  identity from the Ramanujan ``Lost'' Notebook:
$$\prod_{k\geq 1}\frac{1}{(1-q^{5k-4})(1-q^{5k-1})}=\sum_{n\geq 0}
\frac{q^{n^2}}{(1-q)(1-q^2)\dots (1-q^n)} $$
%This identity was claimed to be true by Ramanujan, however
%the proof was found in the earlier publication of Rogers.
Only in the
60s, MacMahon gave a combinatorial interpretation of this identity
as follows.
Let $\lambda=(\lambda_1, \lambda_2, \dots, \lambda_t)$ with
$\sum_i\lambda_i=n$ be a partition of $n$ into non--increasing integers.
Then the identity can derived from the fact that
the number of partitions of $n$ with all $\lambda_i$ of the
 form $5k+1$ or $5k+4$
is equal to the number of partitions where $\lambda_{i}-\lambda_{i+1}\geq 2$
for all $i$.

 This identity is a very special
case of the Andrews identity we used  in \cite{BBL}: For any numbers $b_i,c_i,
i=1,\dots,k$ and positive integer $N$ we have
\begin{multline*}  1+ \sum_{n =1}^N  q^{kn+Nn} (1+ q^n)  \frac{(q^{-N})_n}{(q^{N+1})_n} \prod_{i=1}^k
\frac{(b_i)_{n}}{b_i^{n}} \frac{(c_i)_{n}}{c_i^{n}}
\frac{1}{(\frac{q}{b_i})_n (\frac{q}{c_i})_n}=
\\
\frac{(q)_N\,  (\frac{q}{b_kc_k})_N}{(\frac{q}{b_k})_N\,
(\frac{q}{c_k})_N} \sum_{n_k \ge n_{k-1} \ge \dots \ge n_2 \ge n_1=0}
 \frac{q^{n_k} (q^{-N})_{n_k} (b_k)_{n_k}
(c_k)_{n_k}}{(q^{-N}b_kc_k)_{n_k}}\prod_{i=1}^{k-1}\frac{q^{n_i}
\frac{(b_i)_{n_i}}{b_i^{n_i}} \frac{(c_i)_{n_i}}{c_i^{n_i}}
(\frac{q}{b_i c_i})_{n_{i+1}-n_i} }{(q)_{n_{i+1}-n_i}
(\frac{q}{b_i})_{n_{i+1}}  (\frac{q}{c_i})_{n_{i+1}}} \label{1002}
\end{multline*}
where $(a)_n=(a;q)_n$.

The point is that (for a special choice of parameters)
the left hand side of the identity
 can be identified with
$\cL_{b,c;n}\left((q^{n};q)_{k+1}\; (q^{-n};q)_{k+1}\right)$, where
 the right hand side is a sum with all summands
divisible by $(q^{k+1};q)_{k+1}$.
% but it does  exactly what  we need:
%it shows the required divisibility.

In the case $b=\pm 1$, the computations are especially simple and we get
that
$$\cL_{-1,1;n}((q^{n};q)_{k+1}\; (q^{-n};q)_{k+1})=2(q^{k+1};q)_{k+1}\, .$$
The same holds also for $\cL_{1,1;n}$ up to units.
This allows to write the explicit formula for $I_M$ given in the introduction.

\subsection{Frobenius isomorphism}
It remains to
show that the image of the
Laplace transform
belongs to $\R_b$, i.e. that
certain roots of $q$ exist in $\R_b$.

In \cite{BBL} we proved the following.

\begin{Thm}\label{frob}{\rm [Beliakova--B\"uhler--Le]}
The Frobenius endomorphism
 $F_b: \Z[1/b][q]^{\N_b} \to \Z[1/b][q]^{\N_b}$, sending $q$ to $q^b$,
 is an isomorphism.
\end{Thm}

This implies the
existence of the $b$th root of $q$ in $\calS_{b,0}$
defined by
\[
q^{1/b}:=F^{-1}_b (q)  \in \calS_{b,0}\,.
\]

Let us mention that this result does not hold over $\Q$, i.e.
for  $y \in\Q^{\N_b}$ with $y^b=1$ we have $y=\pm 1$.

Further, we introduce
another Frobenius homomorphism
$$G_m : R[q]^{\N_b} \to R[q]^{m\N_b} \quad \text{ by } \qquad G_m(q) = q^m.$$
Since $\Phi_{mr}(q)$ always divides $\Phi_{r}(q^m)$,  $G_m$ is well--defined.

This map allows us to transfer $q^{1/b}$ from $\calS_{p,0}$
to $\calS_{p,i}$ with $i> 0$, and hence to define a
realization of $q^{a^2/b}$ in $\calS_p$ with the correct evaluation.

\subsection{Integrality}
The proof of integrality of the SO(3) WRT invariants for
all  roots of unity and for all 3--manifolds is given in \cite{BL}.
Note that even restricted to rational homology 3--spheres,
this fact requires a separate proof, since
 the existence
of $I_M \in \R_b$ does not imply  integrality of the WRT invariants,
unless $b= 1$. In this subsection
we define a subring $\Gamma_b \subset \R_b$, such that
for any $f\in \Gamma_b$, $\ev_\xi(f)\in \Z[\xi]$.
%the evaluations
%at roots of unity for the elements of this subring are algebraic integers.
The fact that $I_{M,1}$ belongs
to this subring was proved in \cite{BL}. We conjecture that
this holds in general.

For any divisor $c$ of $b$, let us decompose
$\N=\cup_{c|b}\; c\; \N_{b/c}$. Then
$$\Z[1/b][q]^\N=\prod_{c|b} \Z[1/b][q]^{c\,\N_{b/c}}\, .$$
Analogously, we have
$\Gamma_b=\prod_{c|b} \Gamma_{b,c}$, where
$\Gamma_{b,c} \subset \Z[1/b][q]^{c\, \N_{b/c}}$ is defined as follows.

Let $A_{b,c}=\Z[e_c][q^{\pm 1}, q^{\pm c/b}]$.
Put $t=q^{c/b}$ and let $A^{(m)}_{b,c}$ be the algebra generated over  $A_{b,c}$
by
$$\frac{(t;t)_m}{(q^c;q^c)_{m}}\;
 \frac{(e_c;e_c)_{(c-1)/2}}{\widetilde {(q;q)_m}}\; $$
where $$\widetilde{(q;q)_m}=\prod^m_{i=1, c\, \nmid\,  i}(1-q^i)\, .$$
Then
every element $f\in \Gamma_{b,c}$ has a presentation
$$f=\sum^\infty_{m=0} f_m \; \frac{(q^{m+1};q)_{m+1}}{1-q}
$$
with $f_m\in A^{(m)}_{b,c}$.
For any root of unity $\xi$ of odd order $r$ with $(r,b)=c$
and $f\in \Gamma_{b,c}$, we have
$$ \ev_\xi (f)= \ev_\xi \left(\sum^{(r-3)/2}_{m=0} f_m \frac{(q^{m+1};q)_{m+1}}{1-q}
\right)\; .$$
%where $K$ is the minumum of $(r-1)/2$ and $r/c$.
%(this does not seem very natural, but
% the first divisibility below does not hold if $K=(r-1)/2$,
%and $c$ is big,
%am I right??)

Observe that for $m<(r-1)/2$,
$$
\ev_\xi (\widetilde {(q;q)_{m}}) \; \mid\;
(e_c;e_c)_{(c-1)/2}
\;\;\text{and}\;\;
\ev_\xi \left(\frac{(t;t)_m}{(q^c;q^c)_{m}}
\right) \in \Z[\xi]\, .$$

\begin{conjecture}
For any $M\in \M_b$,
there exists an invariant $I'_M \in \Gamma_b$, such that
for any root of unity $\xi$ of odd order
$\ev_\xi(I'_M)=\tau_M(\xi)$.
\end{conjecture}

We also expect that $\Gamma_{b,b}$ is determined by its Taylor
expansion at $e_b$.

%%%%%%%%%%%%%%%%%%%%%%%%%%%%%%%%%%%%%%%%%%%%%%%%%%%%%%%%%%%%%%%%%%%%%%%%%%%%%%%%%%%%%%%%%%%%%%%%%%%%%%%%%%%%%%%%%%%%%%%%%%%%%%%%%%%%%%%%%%%%%%%%%%%%%%%%%%%%%%%%%%%%%%%%%%%%%%%%%%%%%%%%%%

\end{document}